\documentclass[12pt]{article}
\usepackage{amsmath,amsfonts,amssymb,rotating,colordvi}
\usepackage{colordvi}
 \usepackage{CJK}
\def\qed{\nopagebreak\hfill{\rule{4pt}{7pt}}}
\def\proof{\noindent {\it{Proof.} \hskip 2pt}}

\newtheorem{theo}{Theorem}[section]

\newtheorem{lemm}[theo]{Lemma}

\newtheorem{coro}[theo]{Corollary}

\numberwithin{equation}{section}

\newdimen\Squaresize \Squaresize=11pt
\newdimen\Thickness \Thickness=0.7pt
\def\Square#1{\hbox{\vrule width \Thickness
   \vbox to \Squaresize{\hrule height \Thickness\vss
    \hbox to \Squaresize{\hss#1\hss}
   \vss\hrule height\Thickness}
\unskip\vrule width \Thickness} \kern-\Thickness}

\def\Vsquare#1{\vbox{\Square{$#1$}}\kern-\Thickness}

\def\moins{\raise 1pt\hbox{{$\scriptstyle -$}}}

\begin{document}

\parskip 6pt

\begin{center}
{\large \bf  Ratio Monotonicity of Polynomials

             Derived from   Nondecreasing Sequences}
\end{center}

\begin{center}
William Y. C. Chen$^{1}$, Arthur
L. B. Yang$^{2}$ and Elaine L. F. Zhou$^{3}$\\[6pt]
Center for Combinatorics, LPMC-TJKLC\\
Nankai University, Tianjin 300071, P. R. China

 Email: $^{1}${\tt
chen@nankai.edu.cn}, $^{2}${\tt yang@nankai.edu.cn}, $^{3}${\tt
zhoulf@cfc.nankai.edu.cn}
\end{center}

\vspace{0.3cm} \noindent{\bf Abstract.} The ratio monotonicity of a polynomial is a stronger property than  log-concavity. Let $P(x)$ be a polynomial with nonnegative and nondecreasing coefficients. We prove the ratio monotone property of $P(x+1)$, which leads to the log-concavity of $P(x+c)$ for any $c\geq 1$ due to Llamas and Mart\'{\i}nez-Bernal.  As a consequence, we obtain the ratio monotonicity of the Boros-Moll
polynomials obtained by Chen and Xia without resorting to the recurrence relations of the coefficients.



\noindent {\bf Suggested Running Title:} Ratio Monotonicity

\section{Introduction}

This paper is concerned with the ratio  monotone property of polynomials derived from nonnegative and nondecreasing sequences.
 A sequence $\{a_k\}_{0 \leq k \leq m}$ of positive real
numbers  is said to be  {unimodal} if there exists an integer
$r\geq 0$ such that
$$a_0\leq\cdots\leq a_{r-1}\leq
a_r\geq a_{r+1} \geq \cdots \geq a_m,$$ and it is said to be
 {spiral} if
\begin{equation}
 a_m\leq a_0\leq a_{m-1}\leq a_1 \leq \cdots \leq a_{[\frac{m}{2}]},
\end{equation}
where $[\frac{m}{2}]$ stands for the largest integer less than
$\frac{m}{2}$. We say that a sequence  $\{a_k\}_{0 \leq k \leq m}$ is  {log-concave} if for
any $1\leq k\leq m-1$,
$$a_{k}^2-a_{k+1}a_{k-1}\geq 0,$$
or equivalently,
$$ \frac{a_0}{a_1} \leq \frac{a_1}{a_2} \leq \cdots\leq
 \frac{a_{m-1}}{a_m}.
$$

It is easy to see that either  log-concavity or the spiral property implies
unimodality, while a log-concave sequence is not necessarily spiral,
and vice versa.

A stronger property, which implies both log-concavity and the spiral property,
was introduced by Chen and Xia \cite{chenxia07} and is called the
{ratio monotonicity}. A sequence of positive real numbers
$\{a_k\}_{0 \leq k \leq m} $ is said to be {ratio monotone} if
\begin{equation}
\frac{a_m}{a_0} \leq \frac{a_{m-1}}{a_1} \leq \cdots\leq
 \frac{a_{m-i}}{a_i} \leq \cdots \leq \frac{a_{m-[\frac{m-1}{2}]}}{a_{[\frac{m-1}{2}]}} \leq 1
\end{equation}
and
\begin{equation}
\frac{a_0}{a_{m-1}} \leq \frac{a_{1}}{a_{m-2}} \leq \cdots\leq
 \frac{a_{i-1}}{a_{m-i}} \leq \cdots \leq \frac{a_{[\frac{m}{2}]-1}}{a_{m-[\frac{m}{2}]}} \leq 1.
\end{equation}
Given a
polynomial $ P(x)= a_0 +a_1x+ \cdots + a_mx^m$  with positive coefficients, we
say that $P(x)$ is log-concave (or ratio monotone) if $\{a_k\}_{0
\leq k \leq m}$ is log-concave (resp., ratio monotone).

Assume that $P(x)$ is a polynomial with nonnegative and nondecreasing coefficients.
Boros and Moll \cite{bormol1999}
 proved the unimodality of  $P(x + 1)$ which implies the unimodality of the Boros-Moll polynomials. They posed the conjecture that the Boros-Moll polynomials are log-concave, which was confirmed by  Kauers and Paule \cite{KauPau2007}.
 Alvarez et al. \cite{aabkmr2001}  showed that
$P(x + n)$ is also unimodal for any positive integer $n$. Wang and Yeh \cite{wangye2005}
obtained a stronger result
that $P(x + c)$ is unimodal for $c> 0$. Llamas and Mart\'{\i}nez-Bernal \cite{llamab2010} proved that $P(x + c)$ is log-concave for $c\geq 1$.

In this paper, we prove that if $P(x)$ is a polynomial with nonnegative and nondecreasing coefficients, then $P(x+1)$ is ratio monotone.
This property implies the log-concavity of $P(x+1)$.
Note that by a criterion for log-concavity due to Brenti \cite{brenti1994}, the log-concavity of
$P(x+1)$ leads to the log-concavity of
$P(x+c)$ for $c\geq 1$, as established by Llamas and Mart\'{\i}nez-Bernal \cite{llamab2010}.
The ratio monotonicity of $P(x+1)$ serves as a simple proof
of the ratio monotonicity of the Boros-Moll polynomials obtained by Chen and Xia \cite{chenxia08} without resorting to the recurrence relations of the
coefficients.

\section{The ratio monotone property}

The main result of this paper is given below.

\begin{theo}\label{main}
If $P(x)$ is a polynomial with nonnegative and nondecreasing
coefficients, then $P(x+1)$ is ratio monotone.
\end{theo}

To prove Theorem \ref{main}, we need three lemmas. The first lemma
is a special case of \cite[Lemma 2.1]{chenxia07}.

\begin{lemm}\label{lemm1}
Suppose that $a$, $b$, $c$, $d$, $e$, $f$ are positive real numbers
satisfying
$$\frac{a}{b} \leq \frac{c}{d} \leq \frac{e}{f}.$$
Then
$$\frac{a+c}{b+d} \leq \frac{e+c}{f+d}.$$
\end{lemm}

\begin{lemm}\label{lemm2}
If $B(x)$ is a ratio monotone polynomial, then so is $(x+1)B(x)$.
\end{lemm}

\proof Let \[ B(x)=\sum_{k=0}^ma_kx^k \quad \mbox{ and} \quad  (x+1)B(x)=\sum_{k=0}^{m+1}b_kx^k.\] For each $k$ we have
$b_k=a_{k-1}+a_k$, where $a_{-1}$ and $a_{m+1}$ are set to $0$.

When $m=2n$, the ratio monotonicity of $B(x)$ states that
\begin{equation}\label{e5}
\frac{a_{2n}}{a_0} \leq \frac{a_{2n-1}}{a_1} \leq \cdots \leq
\frac{a_{2n-i}}{a_i} \leq \cdots \leq \frac{a_{n+1}}{a_{n-1}} \leq 1
\end{equation}
and
\begin{equation}\label{e6}
\frac{a_0}{a_{2n-1}} \leq \frac{a_{1}}{a_{2n-2}} \leq \cdots \leq
\frac{a_{i-1}}{a_{2n-i}} \leq \cdots \leq \frac{a_{n-1}}{a_n} \leq
1.
\end {equation}

In order to show that $(x+1)B(x)$ is ratio monotone, we need to
justify that

\begin{equation}\label{e7}
\frac{b_{2n+1}}{b_0} \leq \frac{b_{2n}}{b_1} \leq \cdots \leq
\frac{b_{2n+1-i}}{b_i} \leq \cdots \leq \frac{b_{n+1}}{b_{n}} \leq 1
\end{equation}
and
\begin{equation}\label{e8}
\frac{b_0}{b_{2n}} \leq \frac{b_{1}}{b_{2n-1}} \leq \cdots \leq
\frac{b_{i}}{b_{2n-i}} \leq \cdots \leq \frac{b_{n-1}}{b_{n+1}} \leq
1.
\end {equation}

We first consider \eqref{e7}. Since
$$\frac{a_{2n}}{a_0} \leq
\frac{a_{2n-1}}{a_1},$$ we see that
$$\frac{a_{2n}}{a_0} \leq
\frac{a_{2n-1}+a_{2n}}{a_1+a_0}, $$
that is,
$$ \frac{b_{2n+1}}{b_0}
\leq \frac{b_{2n}}{b_1}.$$ For $1\leq i\leq n-1$, from \eqref{e5} we
deduce that
$$\frac{a_{2n+1-i}}{a_{i-1}}\leq \frac{a_{2n-i}}{a_i}\leq \frac{a_{2n-i-1}}{a_{i+1}}.$$
By Lemma \ref{lemm1}, we obtain
$$\frac{a_{2n+1-i}+a_{2n-i}}{a_i+a_{i-1}}\leq \frac{a_{2n-i}+a_{2n-i-1}}{a_{i+1}+a_i},$$
or equivalently,
 $$\frac{b_{2n+1-i}}{b_i}\leq \frac{b_{2n-i}}{b_{i+1}}.$$
In light of \eqref{e5}, we see that $a_{n+1}\leq a_{n-1}$, and thus we have
$$\frac{b_{n+1}}{b_n}=\frac{a_{n+1}+a_n}{a_n+a_{n-1}}\leq 1.$$

Next, we prove \eqref{e8}. From $\frac{a_0}{a_{2n-1}} \leq
\frac{a_{1}}{a_{2n-2}}$ it follows that
$$\frac{a_0}{a_{2n-1}+a_{2n}}
\leq \frac{a_1+a_0}{a_{2n-2}+a_{2n-1}},$$
that is,
$$\frac{b_0}{b_{2n}} \leq \frac{b_{1}}{b_{2n-1}}.$$
 For $2\leq i\leq n-1$, in view of \eqref{e6} we
find that
$$\frac{a_{i-2}}{a_{2n-i+1}}\leq \frac{a_{i-1}}{a_{2n-i}}\leq \frac{a_{i}}{a_{2n-i-1}}.$$
By Lemma \ref{lemm1}, we have
$$\frac{a_{i-1}+a_{i-2}}{a_{2n-i+1}+a_{2n-i}}\leq \frac{a_{i}+a_{i-1}}{a_{2n-i}+a_{2n-i-1}},$$
which can be expressed as
$$\frac{b_{i-1}}{b_{2n-i+1}}\leq \frac{b_{i}}{b_{2n-i}}.$$
 From \eqref{e6} it is clear that $a_{n-2}\leq a_{n+1}$ and $a_{n-1}\leq a_n$, and hence
$$\frac{b_{n-1}}{b_{n+1}}=\frac{a_{n-1}+a_{n-2}}{a_{n+1}+a_n}\leq 1.$$

The case $m=2n+1$ can be dealt with in the same manner. This completes the proof.
 \qed

The third lemma is concerned with an inequality of an increasing
positive sequence. This inequality will be needed in the proof of Theorem \ref{main}.

\begin{lemm}\label{lemm3}
For any nondecreasing positive sequence $\{a_{k}\}_{0\leq k\leq m}$,
we have
$$\frac{m(m+1)}{2}a_m^2+a_ma_{m-1}\geq \left(\sum_{k=0}^{m-2}
\left(m-1-k\right)a_k\right)a_{m-1}+\left(\sum_{k=0}^ma_k\right)a_{m-2}.$$
\end{lemm}
\proof Since $0 < a_0\leq a_1\leq  \cdots \leq a_{m-1} \leq
 a_m,$ we have
\begin{eqnarray*}
\lefteqn{\frac{m(m+1)}{2}a_m^2+a_ma_{m-1}-\left(\sum_{k=0}^{m-2}
\left(m-1-k\right)a_k\right)a_{m-1}-\left(\sum_{k=0}^ma_k\right)a_{m-2}}\\
&\geq & \frac{m(m+1)}{2}a_m^2+a_ma_{m-1}-\sum_{k=0}^{m-2}
(m-1-k)a_m^2-\sum_{k=1}^m a_m^2-a_ma_{m-1} ,
\end{eqnarray*}
which simplifies to zero, as desired. \qed

\noindent \textit{Proof of Theorem \ref{main}.} Use induction on the degree
$m$ of $P(x)$. Let
$$P(x)= \sum_{k=0}^{m}a_kx^k,$$
where $0 < a_0\leq a_1\leq  \cdots \leq a_{m-1} \leq
 a_m.$

When $m=2$, we have
$$P(x+1)=a_2x^2+(a_1+2a_2)x +a_0+a_1+a_2.$$
Note that
$
a_2 \leq a_0+a_1+a_2$, $a_0+a_1+a_2 \leq a_1+2a_2.
$
Therefore, the theorem holds for $m=2$.

Now assume that the theorem holds for polynomials of degree $m-1$.
We need to show that it is also true for polynomials $P(x)$ of
degree $m$. Suppose  that
\begin{equation}\label{b0exp}
P(x+1) =\sum_{k=0}^{m}a_k(z+1)^k= \sum_{k=0}^{m}b_kz^k.
\end{equation}
We wish to prove that
\begin{equation}\label{e9}
\frac{b_m}{b_0} \leq \frac{b_{m-1}}{b_1} \leq \cdots\leq
 \frac{b_{m-i}}{b_i} \leq \cdots \leq \frac{b_{m-[\frac{m-1}{2}]}}{b_{[\frac{m-1}{2}]}} \leq 1
\end{equation}
and
\begin{equation}\label{e10}
\frac{b_0}{b_{m-1}} \leq \frac{b_{1}}{b_{m-2}} \leq \cdots \leq
 \frac{b_{i-1}}{b_{m-i}} \leq \cdots \leq \frac{b_{[\frac{m}{2}]-1}}{b_{m-[\frac{m}{2}]}} \leq 1.
\end{equation}

Let $$Q(x)= \sum_{k=0}^{m-1}a_{k+1}x^k.$$ Then
$$P(x+1) = a_0+(x+1)Q(x+1).$$
By the inductive hypothesis and Lemma \ref{lemm2}, we deduce that the polynomial
$$(x+1)Q(x+1)=b_0-a_0+\sum_{k=1}^mb_kx^k$$
is ratio monotone. So we have
\begin{equation}\label{e11}
\frac{b_m}{b_0-a_0} \leq \frac{b_{m-1}}{b_1} \leq \cdots\leq
 \frac{b_{m-i}}{b_i} \leq \cdots \leq \frac{b_{m-[\frac{m-1}{2}]}}{b_{[\frac{m-1}{2}]}} \leq 1
\end{equation}
and
\begin{equation}\label{e12}
\frac{b_0-a_0}{b_{m-1}} \leq \frac{b_{1}}{b_{m-2}} \leq \cdots\leq
 \frac{b_{i-1}}{b_{m-i}} \leq \cdots \leq \frac{b_{[\frac{m}{2}]-1}}{b_{m-[\frac{m}{2}]}} \leq 1.
\end{equation}

Clearly, \eqref{e9} follows from \eqref{e11}. To prove
\eqref{e10}, it remains to show that
$$\frac{b_0}{b_{m-1}} \leq \frac{b_{1}}{b_{m-2}}.$$

From \eqref{b0exp}, we see that
\begin{align*}
b_0 & =  \sum_{k=0}^ma_k, \qquad b_{m-1}=a_{m-1}+ma_m,
\end{align*}
and
\begin{align*}
b_1 & =  \sum_{k=0}^m ka_k, \quad \,\,\,b_{m-2} =
a_{m-2}+(m-1)a_{m-1}+\binom{m}{2}a_m.
\end{align*}

Consequently, it suffices to show that
$$ \frac{\sum_{k=0}^ma_k}{a_{m-1}+ma_m} \leq \frac{\sum_{k=0}^m
ka_k}{a_{m-2}+(m-1)a_{m-1}+\binom{m}{2}a_m},
$$
or equivalently,
\begin{eqnarray*}
&&\left(\sum_{k=0}^m ka_k\right)a_{m-1}+\left(\sum_{k=0}^m
mka_k\right)a_m-\left(\sum_{k=0}^ma_k\right)a_{m-2}\\
&&\quad -\left(\sum_{k=0}^m(m-1)a_k\right)a_{m-1}
-\left(\sum_{k=0}^m\binom{m}{2}a_k\right)a_m \geq 0.
\end{eqnarray*}
The left hand side of the above inequality can be simplified to
$$
\left(\sum_{k=0}^m\frac{2k-m+1}{2}a_k\right)ma_m+
\left(\sum_{k=0}^m
\left(k-m+1\right)a_k\right)a_{m-1}-\left(\sum_{k=0}^ma_k\right)a_{m-2},
$$
which can be rewritten as a sum of
\begin{equation}
\left(\sum_{k=0}^{m-1}\frac{2k-m+1}{2}a_k\right)ma_m \label{s1}
\end{equation}
and
\begin{align}
\frac{m(m+1)}{2}a_m^2+a_ma_{m-1}-\left(\sum_{k=0}^{m-2}
\left(m-1-k\right)a_k\right)a_{m-1}-\left(\sum_{k=0}^ma_k\right)a_{m-2}.\label{s2}
\end{align}
%
By Lemma \ref{lemm3}, the sum in \eqref{s2} is nonnegative. The sum in \eqref{s1} is also nonnegative, since
\begin{eqnarray*}
\sum_{k=0}^{m-1}\frac{2k-m+1}{2}a_k&=&\sum_{k=[\frac{m-1}{2}]+1}^{m-1}\frac{2k-m+1}{2}a_k
-\sum_{k=0}^{[\frac{m-1}{2}]}\frac{m-1-2k}{2}a_k\\
&=&\sum_{k=0}^{m-2-[\frac{m-1}{2}]}\frac{m-1-2k}{2}a_{m-1-k}
-\sum_{k=0}^{[\frac{m-1}{2}]}\frac{m-1-2k}{2}a_k\\
&=&\sum_{k=0}^{[\frac{m-1}{2}]}\frac{m-1-2k}{2}(a_{m-1-k}-a_k),
\end{eqnarray*}
which is nonnegative, and thus the proof is complete. \qed

 Theorem \ref{main} leads to the following  result of Llamas and Mart\'{\i}nez-Bernal \cite{llamab2010}, since the ratio monotonicity
 implies log-concavity of $P(x+1)$ and the log-concavity of $P(x+1)$ implies the log-concavity of $P(x+c)$ for $c\geq 1$ by a
 criterion of Brenti \cite{brenti1989, brenti1994}.

\begin{coro}
If $P(x)$ is a polynomial with nonnegative and nondecreasing
coefficients, then for any $c\geq 1$ the polynomial $P(x + c)$ is
log-concave and has no internal zero coefficients.
\end{coro}

Theorem \ref{main} also serves as a simple proof of the ratio monotonicity of the Boros-Moll polynomials $P_m(x)$, which were introduced by Boros and
Moll \cite{borosmoll1999} in their study of the following
quartic integral
$$
\int_{0}^{+\infty}\frac{1}{(t^4+2xt^2+1)^{m+1}}dt=\frac{\pi}{2^{m+3/2}(x+1)^{m+1/2}}P_m(x).
$$
Let
$$c_k(m)=2^{-2m+k}\binom{2m-2k}{m-k}\binom{m+k}{k}.$$
Boros and Moll showed that
\begin{equation}\label{bmpol}
P_m(x)= \sum_{k=0}^m c_k(m) (x+1)^k.
\end{equation}
They also observed that, for $0\leq k\leq m-1$,
$$\frac{c_{k}(m)}{c_{k+1}(m)}=\frac{(2m-2k-1)(k+1)}{(m-k)(m+k+1)}<1.$$
Thus, $P_m(x-1)$ is a polynomial with nonnegative and nondecreasing coefficients. Boros and Moll \cite{borosmoll1999}
proved that $P_m(x)$ is unimodal for any $m\geq¡¡0$, and Moll \cite{moll2002} conjectured that $P_m(x)$ is log-concave for any $m$. This conjecture was confirmed by Kauers and Paule \cite{KauPau2007}. The ratio monotonicity of $P_m(x)$ was established by Chen and Xia and the proof is quite technical and heavily depends on inequalities on the coefficients. The proof of Theorem \ref{main} shows that
the log-concavity and ratio monotonicity only depend on the nondecreasing property of the coefficients of $P_m(x-1)$.

\vskip 3mm

\noindent {\bf Acknowledgments.} This work was supported by the 973
Project, the PCSIRT Project of the Ministry of Education, and the
National Science Foundation of China.


\begin{thebibliography}{99}

\bibitem{aabkmr2001} J. Alvarez, M. Amadis, G. Boros, D. Karp, V.H. Moll and L. Rosales,
{An extension of a criterion for unimodality}, Electron. J.
Combin. 8 (2001), \#R30.

\bibitem{borosmoll1999} G. Boros and V.H. Moll,
{A sequence of unimodal polynomials}, J. Math. Anal. Appl.
237 (1999), 272--285.


\bibitem{bormol1999} G. Boros and V.H. Moll.
{A criterion for unimodality}, Electron. J. Combin. 6 (1999),
\#R10.

\bibitem{brenti1989} F. Brenti, {Unimodal, log-concave, and P$\mathrm{\acute{o}}$lya frequency sequences in
combinatorics}, Mem. Amer. Math. Soc. 413 (1989), 1--106.


\bibitem{brenti1994} F. Brenti, {Log-concave and unimodal sequences in algebra, combinatorics and geometry: an update},
Contemp. Math. {178} (1994), 71--89.

\bibitem{chenxia07}W.Y.C. Chen and E.X.W. Xia,
{The ratio monotonicity of the $q$-derangement numbers},
arXiv:math.CO/0708.2532.

\bibitem{chenxia08}W.Y.C. Chen and E.X.W. Xia,
{The ratio monotonicity of the Boros-Moll polynomials},
Math. Comput. 78 (2009), 2269--2282.

\bibitem{KauPau2007} M. Kauers and P. Paule, {A computer proof of Moll's log-concavity conjecture},
Proc. Amer. Math. Soc. 135 (2007), 3847--3856.

\bibitem{llamab2010} A. Llamas, J. Mart\'{\i}nez-Bernal, {Nested log-concavity}, Commun. Algebra 38 (2010), 1968--1981.

\bibitem{moll2002} V.H. Moll, {The evaluation of integrals: A personal story}, Notices
Amer. Math. Soc. 49 (2002), 311--317.

\bibitem{stanley1989} R.P. Stanley, {Log-concave and unimodal sequences in algebra, combinatorics and
geometry}, Ann. New York Acad. Sci {576} (1989), 500--535.

\bibitem{wangye2005} Y. Wang and Y.-N. Yeh, {Proof of a conjecture on unimodality},
European J. Combin. 26 (2005), 617--627.

\end{thebibliography}
\end{document}